\documentclass[12pt]{article}
\usepackage[cm]{fullpage}

\usepackage{amsmath,amscd, float,rotating}
\usepackage{pb-diagram}

\usepackage{latexsym}
\usepackage{amsfonts}
\usepackage{amsmath}
\usepackage{amssymb}

\numberwithin{equation}{section}
\newcommand{\qed}{\hfill \ensuremath{\Box}}

\def\XXint#1#2#3{{\setbox0=\hbox{$#1{#2#3}{\int}$}
\vcenter{\hbox{$#2#3$}}\kern-.5\wd0}}

\newcommand{\kod}{\textnormal{kod}}

\newcommand{\ga}{\gamma}
\newcommand{\de}{\delta}

\newcommand{\dbar}{\overline{\partial}}

\newcommand{\ddt}[1]{\frac{\partial #1}{\partial t}}

\newcommand{\ddbar}{\sqrt{-1}\partial\dbar}

\begin{document}
\newcounter{remark}
\newcounter{theor}
\setcounter{remark}{0} \setcounter{theor}{1}
\newtheorem{claim}{Claim}
\newtheorem{theorem}{Theorem}[section]
\newtheorem{proposition}{Proposition}[section]
\newtheorem{lemma}{Lemma}[section]
\newtheorem{conjecture}{Conjecture}[section]
\newtheorem{definition}{Definition}[section]
\newtheorem{corollary}{Corollary}[section]
\newenvironment{proof}[1][Proof]{\begin{trivlist}
\item[\hskip \labelsep {\bfseries #1}]}{\end{trivlist}}
\newenvironment{remark}[1][Remark]{\addtocounter{remark}{1} \begin{trivlist}
\item[\hskip \labelsep {\bfseries #1
\thesection.\theremark}]}{\end{trivlist}}

\begin{center}
{\bf\Large Finite time extinction of the K\"ahler-Ricci flow \footnote{Research supported in part by National Science Foundation
grant  DMS-06-04805.}   }

\bigskip
\bigskip

{\large Jian Song$^{*}$ } \\

\end{center}

\bigskip

\bigskip
\bigskip
\noindent
{\bf Abstract} \ We investigate the limiting behavior of the unnormalized K\"ahler-Ricci flow  on a K\"ahler manifold with a polarized initial K\"ahler metric.  We prove that the K\"ahler-Ricci flow becomes extinct in finite time if and only if the manifold has positive first Chern class and the initial K\"ahler class is proportional to the first Chern class of the manifold.   This proves a conjecture of Tian for the smooth solutions of the K\"ahler-Ricci flow.

\section{ Introduction}

The limiting behavior of the K\"ahler-Ricci flow is deeply related to the existence of canonical metrics and the minimal model program in algebraic geometry.
Let $X$ be an $n$-dimensional compact K\"ahler manifold. We consider the following unnormalized K\"ahler-Ricci flow starting with a K\"ahler metric $\omega_0$.
\begin{equation}\label{ricciflow}
\left\{
\begin{array}{rl}
&{ \displaystyle \ddt{}\omega =-Ric(\omega) }\\
&\\
& \omega|_{t=0}= \omega_0 .
\end{array} \right.
\end{equation}

If  $c_1(X)=0$, the unnormalized K\"ahler-Ricci flow (\ref{ricciflow}) converges to the unique Ricci-flat metric in $[\omega_0]$ \cite{Cao}, \cite{Y1}. If $c_1(X)<0$, the normalized K\"ahler-Ricci flow
\begin{equation} \label{hyperflow}
\ddt{}\omega =-Ric(\omega)  - \omega
\end{equation}
converges to the unique hyperbolic K\"ahler-Einstein metric \cite{Cao}, \cite{Y1}, \cite{A}.

If $c_1(X)$ is nonpositive  and not strictly negative, there does not exist a smooth K\"ahler-Einstein metric on $X$. However, the normalized K\"ahler-Ricci flow (\ref{hyperflow}) has long time existence. If $X$ is a minimal model of general type, the K\"ahler-Ricci flow (\ref{hyperflow})  converges to the unique singular K\"ahler-Einstein metric as time tends to infinity \cite{Ts}, \cite{TZha}. In general, if the canonical bundle $K_X$ is semi-ample, the normalized K\"ahler-Ricci flow (\ref{hyperflow}) can possibly collapse but it always converges to the unique generalized K\"ahler-Einstein metric on the canonical model of $X$ \cite{SoT1}, \cite{SoT2}.

If $c_1(X)>0$, $X$ is called a Fano manifold in algebraic geometry.  It is conjectured by Yau \cite{Y2} that the existence of a K\"ahler-Einstein metric is equivalent to suitable stability in the sense of geometric invariant theory. The condition of $K$-stability is later proposed by Tian \cite{T1} and has been refined  by Donaldson \cite{D}. The K\"ahler-Ricci flow provides an approach to the Yau-Tian-Donaldson conjecture  and it has attracted considerable current interest. We refer the reader to an incomplete list of literatures \cite{PS1}, \cite{PS2}, \cite{PSSW1}, \cite{PSSW2}, \cite{Sz},  \cite{To} and \cite{CW} for some recent development.  If one assumes the existence of  a K\"ahler-Einstein metric, then according to the unpublished work of Perelman \cite{P2} (see \cite{TZhu}), the normalized K\"ahler-Ricci flow $\ddt{}g= -Ric(g) +g$ with the initial K\"ahler metric in $ c_1(X)$  converges to a K\"ahler-Einstein metric (this is due to \cite{H}, \cite{Ch} in the case of one complex dimension).

Let $M$ be a closed Riemannian manifold and $g(t)$ be a smooth solution of the Ricci flow for $t\in [0, T)$. The Ricci flow becomes extinct at $t=T$ if $(M, g(t))$ converges to a point in the sense of Gromov-Hausdorff, or equivalently, the diameter of $g(t)$ tends to $0$ as $t\rightarrow T$. It is proved in \cite{IK} that  the diameter of $g(t)$ is uniformly bounded below from $0$ if there exists an element of infinite order in $H_1(M, \mathbf{Z})$. However, many projective manifolds do not admit elements of infinite order in $H_1(M, \mathbf{Z})$. In particular, Fano manifolds are  all simply connected. Indeed, it is proved by Perelman \cite{P2} (see \cite{SeT}) that the K\"ahler-Ricci flow (\ref{ricciflow}) on Fano manifolds becomes extinct in finite time if the initial K\"ahler class is proportional to the first Chern class of the manifold. Thus it is of our interest to completely understand  the extinction conditions for the K\"ahler-Ricci flow (\ref{ricciflow}). In general, the K\"ahler-Ricci flow can be weakly defined on singular projective varieties with mild singularities and surgeries might be performed, as investigated in \cite{SoT3}. Therefore the notion of extinction can be generalized to the K\"ahler-Ricci flow with surgery. The following conjecture is proposed by Tian in \cite{T2}.

\begin{conjecture} \label{conj1}
The K\"ahler-Ricci flow (\ref{ricciflow}) with surgery becomes extinct in finite time if
and only if the initial projective manifold is birationally equivalent to a Fano manifold.

\end{conjecture}

It is shown in \cite{SW} that the conjecture holds in the special case of Hirzebruch surfaces when the initial K\"ahler metric satisfies the Calabi symmetry.

We prove Conjecture \ref{conj1} in the case of the usual K\"ahler-Ricci flow with smooth solutions.

\begin{theorem}\label{mainth1} Let $X$ be an $n$-dimensional K\"ahler manifold. Then the K\"ahler-Ricci flow (\ref{ricciflow}) with an initial K\"ahler metric in the class of $H^2(X, \mathbf{Z})$ becomes extinct in finite time if and only if $X$ is Fano and the initial K\"ahler class is proportional to $c_1(X)$.

\end{theorem}

More generally, we have the following general diameter estimates.

\begin{theorem} \label{mainth2} Let $X$ be $n$-dimensional K\"ahler manifold. Let $\omega_0$ be a K\"ahler metric on $X$ such that $[\omega_0]\in H^2(X, \mathbf{Z})$. Then the  diameter  is uniformly bounded below from $0$ along the K\"ahler-Ricci flow (\ref{ricciflow}) if one of the following conditions holds.

\begin{enumerate}

\item  $K_X$ is semi-ample.

\item $K_X$ is not nef and  $K_X^{-1}$ is not ample.

\item $K_X^{-1}$ is ample and $[\omega_0]$ is not proportional to $[K^{-1}_X]$.

\end{enumerate}

In particular, if $K_X$ is semi-ample and the Kodaira dimension of $X$ is positive, the diameter tends to infinity of order $\sqrt{t}$ as $t\rightarrow \infty$ along the K\"ahler-Ricci flow (\ref{ricciflow}).

\end{theorem}

It is known that the abundance conjecture holds for dimension three and so $K_X$ is semi-ample whenever it is nef. The following corollary is then an immediate consequence of Theorem \ref{mainth2}.
\begin{corollary} \label{mainth3} Let $X$ be a K\"ahler manifold of $\dim X\leq 3$.  Then the K\"ahler-Ricci flow (\ref{ricciflow}) on $X$ with an initial K\"ahler metric in $H^2(X, \mathbf{Z})$ becomes extinct if and only if $X$ is Fano and the initial K\"ahler class is proportional to $c_1(X)$.

\end{corollary}

The above statement holds for $\dim X\geq 4$ if we assume the abundance conjecture.

We further propose the following conjecture as a natural attempt to generalize Theorem \ref{mainth1} and Corollary \ref{mainth3} by removing the condition for the initial K\"ahler class and the assumption of the abundance conjecture.

\begin{conjecture} Let $X$ be an $n$-dimensional K\"ahler manifold. Then the K\"ahler-Ricci flow on $X$ with an initial K\"ahler metric becomes extinct if and only if $X$ is Fano and the initial K\"ahler class is proportional to $c_1(X)$.

\end{conjecture}

\begin{remark} Theorem 1.1 also holds if the initial K\"ahler metric is in a multiple of a class in $H^2(X, \mathbf{Z})$. This can be easily proved by suitable scaling of the unnormalized flow (\ref{ricciflow}). It is pointed out by V. Tosatti to the author that Theorem 1.1 should also be true if the initial K\"ahler metric sits in the real N\'eron-Severi group $(H^2(X, \mathbf{Z})\cap H^{1,1}(X, \mathbf{C}))\otimes \mathbf{R}$ as the base-point-free theorem holds for $\mathbf{R}$-divisors due to Shokurov.

\end{remark}

\section{Base point freeness}

Let $X$ be an $n$-dimensional projective manifold and $ L \rightarrow X $ a holomorphic line bundle over $X$. Let $N( L) $ be the semi-group defined by $$N(L) = \{ m\in \mathbf{N} ~ | ~ H^0(X , L^m) \neq  0 \}.$$

Given any $m\in N(L)$,  the linear system $|L^m|= \mathbf{P} H^0(X, L^m)$ induces a rational map
$$\Phi_m ~ : ~ X  \dashrightarrow \mathbf{CP}^{d_m}$$
by any basis $\{ \sigma_{m,0},~ \sigma_{m,1}, ~... ~, \sigma_{m, d_m} \}$ of  $H^0(X, L^m )$ as

$$\Phi_m (z) = \left[ \sigma_{m,0},~ \sigma_{m,1}, ~... ~, \sigma_{m, d_m}\right]  (z) , $$
where $d_m + 1 = \dim H^0(X, L^m) $.
Let $Y_m = \overline{\Phi_m (X)} \subset \mathbf{CP}^{d_m}$ be the closure of the image of $X$ by $\Phi_m$.

\begin{definition} The Iitaka dimension of $L$ is defined to be

$$\kappa(X, L) = \max _{ m\in N(L)} \{ \dim Y_m \}$$ if $N(L) \neq \phi$,  and
$\kappa(X, L) = -\infty $ if $N(L)= \phi$.

\end{definition}

\begin{definition}
Let $X$ be an projective manifold and $K_X$ the canonical line bundle over $X$. Then the Kodaira dimension $\kod(X)$ of $X$ is defined to be

$$\kod(X) = \kappa (X, K_X).$$

\end{definition}

The Kodaira dimension  is a birational invariant of projective varieties and the Kodaira dimension of a singular variety is  equal to that of its smooth model.

\begin{definition}

Let $L \rightarrow X$ be a holomorphic line bundle  over a projective manifold $X$. $L$ is called semi-ample if the linear system $| L^m |$ is base point free for some $m>0$.

\end{definition}

For any $m\in \mathbf{N}$ such that $|L^m|$ is base point free,  the linear system $|L^m|$ induces a holomorphic map $\Phi_m$
$$\Phi_m ~ : ~ X  \rightarrow \mathbf{CP}^{d_m}$$
by any basis  of  $H^0(X, L^m )$.
Let $Y_m = \Phi_m (X)$ and so

$$ \Phi_m ~ : ~ X  \rightarrow Y_m \in \mathbf{CP}^{d_m}.$$

The following theorem is well-known (see \cite{L}, \cite{U}).

\begin{theorem}\label{safibration}

Let $L \rightarrow X$ be a semi-ample line bundle over an algebraic manifold $X$. Then there is a projective fibration
$$\pi : X \rightarrow Y $$
such that for any sufficiently large integer $ m $ with $L^m$ being globally generated,
$$Y_m = Y ~~~~ and ~~~~ \Phi_m = \pi, $$
where $Y$ is a normal projective variety.
Furthermore, there exists an ample line bundle $A$ on $Y$ such that $ L^m = \pi^* A$.

\end{theorem}

If $L$ is semi-ample, the graded ring $R(X, L) = \oplus_{m\geq 0} H^0( X, L^m)$ is finitely generated and so $R(X , L)=\oplus_{m\geq 0} H^0(X, L^m)$ is the coordinate ring of $Y$.

\begin{definition} \label{iitaka}
Let $L \rightarrow X$ be a semi-ample line bundle over a projective manifold $X$. Then the algebraic fibration
$\pi: X \rightarrow Y$  as in Theorem \ref{safibration} is called the Iitaka fibration associated to $L$. It is completely determined by the linear system $|L^m|$ for sufficiently large $m$.

\end{definition}


The following theorems are known as the rationality theorem and base-point-free theorem in the minimal model program (see \cite{KMM}, \cite{KM}).
\begin{theorem}\label{rationality} Let $X$ be a projective manifold such that $K_X$ is not nef. Let $H$ be an ample divisor and let

\begin{equation}
\lambda = \max \{ t\in \mathbf{R}~|~ H + t K_X ~ is ~ nef~\}.
\end{equation}
Then $\lambda\in \mathbf{Q}$.

\end{theorem}


\begin{theorem}\label{basepointfree} Let $X$ be a projective manifold. Let $D$ be a nef divisor such that $aD - K_X$ is nef and big for some $a>0$. Then $D$ is semi-ample.

\end{theorem}


We now will apply the base-point-free theorem to the K\"ahler-Ricci flow at the finite blow-up time. We consider the unnormalized K\"ahler-Ricci flow (\ref{ricciflow}) with  the initial K\"ahler class $H = [\omega_0 ] \in H^2(X, \mathbf{Z})\cap H^{1,1}(X, \mathbf{C})$.
The evolution of the K\"ahler class satisfies the following ordinary differential equation

\begin{equation}
\ddt{} [ \omega] = [K_X],     ~~~~~~~~[\omega] |_{t=0} = [H].
\end{equation}
Therefore, $ [ \omega] = [H+ t K_X]$ for $t\geq 0$.


\begin{lemma}

Let
\begin{equation}
T = \sup \{ t \geq 0 ~|~ H + tK_X ~is ~nef\}.
\end{equation}
Then  $T \in \mathbf{Q}$ if $K_X$ is not nef and  $T = \infty$ if $K_X$ is nef.
\end{lemma}

\begin{proof} The lemma is an immediate corollary of the rationality Theorem \ref{rationality}.

\qed
\end{proof}


The following theorem is proved in \cite{TZha} for the maximal existence of the K\"ahler-Ricci flow.
\begin{proposition} \label{maxtime}The K\"ahler-Ricci flow (\ref{ricciflow}) exists for $t\in [0, T)$.

\end{proposition}


\begin{lemma} Let $ L = H + T K_X$. Then $L$ is semi-ample.

\end{lemma}

\begin{proof} Notice that $L - \epsilon K_X = H + (T-\epsilon) K_X$ is ample for sufficiently small $\epsilon>0$. The lemma follows from the base-point-free theorem \ref{basepointfree}.

\qed
\end{proof}


\section{The case when $K_X$ is not nef and $\kappa(L)=0$}

If $K_X$ is not nef, $T< \infty$ and the unnormalized K\"ahler-Ricci flow (\ref{ricciflow}) does not have long time existence and it must develop singularities at $t=T$. In particular, it is shown in \cite{Z} that the scalar curvature must blow up at $t=T$. Let $L= H+T K_X$. Then either $\kappa(L)=0$ or $\kappa(L)>0$ as $L$ is semi-ample.

\begin{proposition}\label{fano} If $\kappa(L)=0$, $X$ is Fano and $ [H] $ is proportional to $c_1(X)$.

\end{proposition}

\begin{proof} Since $L$ semi-ample and $\kappa(L)=0$, $L^m$ is trivial for some $m$ and $[L]=0$. Hence $$T ~c_1(X)= -T ~[K_X] = [H]>0.$$

\qed\end{proof}

The following theorem is proved by Perelman \cite{P2} (see \cite{SeT}).

\begin{theorem} \label{perelman}Let $X$ be a Fano manifold of complex dimension $n\geq 2$. The unnormalized K\"ahler-Ricci flow becomes extinct in finite time if the initial K\"ahler class is proportional to the first Chern class $c_1(X)$.

\end{theorem}

More precisely, Perelman shows that if $X$ is Fano and the initial K\"ahler metric lies in $c_1(X)$, the diameter of the evolving metrics is uniformly bounded above along the normalized K\"ahler-Ricci flow $\ddt{}g= -Ric(g) +g$. The above theorem follows immediately after scaling the normalized flow back to the unnormalized flow.


\section{The case when $K_X$ is not nef and $\kappa(L)>0$}

Now we assume $ k = \kappa(L) >0$. Since $L$ is semi-ample, $H^0(X, L^m)$ induces a holomorphic map for sufficiently large $m$ $$\pi: X \rightarrow Y  \subset \mathbf{CP}^{d_m}, $$
where $Y$ is a normal variety of $\dim Y = k$ and $d_m+1 = \dim H^0(X, L^m)$.

Since $L=H+ T K_X$ is semi-ample, there exists a smooth $(1,1)$-form $\chi\in [K_X]$ such that
$$\omega_{T} = \omega_0 + T \chi \geq 0$$ is proportional to the pullback of the Fubini-Study metric $\omega_{FS}$ on $\mathbf{CP}^{d_m}$ by $\Phi$.
There also exists a smooth volume form $\Omega$ on $X$ such that $$\ddbar \log \Omega = \chi. $$

Let $ \omega_t = \omega_0 + t \chi$. Then

$$ \omega_t = \frac{t}{T} \omega_T + \frac{T-t}{T} \omega_0\geq \frac{T-t}{T} \omega_0 >0$$
and the following lemma holds immediately.
\begin{lemma} For any $(t, z) \in [0, T) \times X$,

\begin{equation}
\frac{ \omega_t^n} { \omega_0^n} \geq (\frac{T-t}{T})^n.
\end{equation}

\end{lemma}


We consider the following Monge-Amp\`ere flow induced by the unnormalized K\"ahler-Ricci flow.

\begin{equation}\label{maflow}
\left\{
\begin{array}{rl}
&{ \displaystyle \ddt{}\varphi = \log \frac{  ( \omega_t + \ddbar \varphi)^n } { \Omega}   } \\
&\\
& \varphi |_{t=0}= 0 .
\end{array} \right.
\end{equation}

\begin{lemma}

There exists $C>0$ such that for $ (t, z)\in [0, T) \times X$,

\begin{equation}
\varphi\leq C.
\end{equation}

\end{lemma}

\begin{proof} Let $\varphi_{max} (t) = \max_{z\in X} \varphi(t, z)$. Then

$$\ddt{} \varphi_{max} \leq \log \frac{\omega_t^n}{\Omega} \leq \log \frac{ (\omega_0 + \omega_T)^n}{\Omega } \leq C.$$

Hence $$ \varphi(t, z)  \leq  CT.$$

\qed
\end{proof}


\begin{proposition}   There exists $C>0$ such that for all $(t, z) \in [0, T) \times X  $,

\begin{equation}
tr_{\omega} (\omega_T) \leq C.
\end{equation}

\end{proposition}

\begin{proof}  This is a parabolic Schwarz lemma similar to the one given in \cite{SoT1}. Suppose $\omega_T= c \pi^* \omega_{FS}$ and $\omega_{FS} = h_{\alpha\bar{\beta}} d x^\alpha\wedge d\overline{x}^\beta$.

Choose normal coordinate systems for $g=\omega(t,\cdot)$ on $X$
and $h$ on $\mathbf{CP}^{d_m}$ respectively. Let
$u=tr_{g}(h)=g^{i\overline{j}} \pi ^{\alpha}_i
\pi^{\overline{\beta}}_{\overline{j}}h_{\alpha\overline{\beta}}$ and
we will calculate the evolution of $u$. $u$ is nonnegative as $\pi$ is holomorphic. Standard calculation shows that
\begin{eqnarray*}\label{sch}
\Delta u &=& g^{k\overline{l}}\partial_k
\partial_{\overline{l}} \left( g^{i\overline{j}}\pi^{\alpha}_i
\pi^{\overline{\beta}}_{\overline{j}}h_{\alpha\overline{\beta}} \right)\\
&=&g^{i\overline{l}}g^{k\overline{j}}R_{k\overline{l}}
\pi^{\alpha}_{i}
\pi^{\overline{\beta}}_{\overline{j}}h_{\alpha\overline{\beta}}+
g^{i\overline{j}}g^{k\overline{l}}\pi^{\alpha}_{i,k}\pi^{\overline{\beta}}_{\overline{j},\overline{l}}h_{\alpha\overline{\beta}}
-g^{i\overline{j}}g^{k\overline{l}}S_{\alpha\overline{\beta}
\ga\overline{\de}} \pi^{\alpha}_i
\pi^{\overline{\beta}}_{\overline{j}}\pi^{\ga}_k
\pi^{\overline{\de}}_{\overline{l}},
\end{eqnarray*}
where $S_{\alpha\overline{\beta} \ga\overline{\de}}$ is the
curvature tensor of $h_{\alpha\bar{\beta}}$.
 \noindent By the definition of $u$ we have
\begin{eqnarray*}
\Delta u \geq g^{i\overline{l}}g^{k\overline{j}}R_{k\overline{l}}
\pi^{\alpha}_i
\pi^{\overline{\beta}}_{\overline{j}}h_{\alpha\overline{\beta}}- Ku^2
\end{eqnarray*}
for some fixed constant $K>0$. Now
\begin{eqnarray*}
\ddt{u}&=&-g^{i\overline{l}}g^{k\overline{j}}\ddt{g_{k\overline{l}}}
\pi^{\alpha}_i
\pi^{\overline{\beta}}_{\overline{j}} h_{\alpha\overline{\beta}}\\
&=&g^{i\overline{l}}g^{k\overline{j}}R_{k\overline{l}}\pi^{\alpha}_i
\pi^{\overline{\beta}}_{\overline{j}}h_{\alpha\overline{\beta}},
\end{eqnarray*}
therefore
\begin{equation}
\left( \ddt{}-\Delta \right) \log  u\leq Ku.
\end{equation}

Consider $H = \log u - 2A\varphi$. If $A$ is chosen to be sufficiently large,

\begin{eqnarray*}
&&(\ddt{} - \Delta) H \\
&\leq& - tr_{\omega} (2A\omega_t - K \omega_T) -2 A \log \frac{\omega^n}{\Omega} + 2An\\
&\leq&  - tr_{\omega} (A\omega_t) + 2A\log  \frac{\omega_t^n}{\omega^n}  + 2A\log \frac{\Omega}{\omega_t^n} + 2An\\
&\leq & - CA n \log (T-t) + CA .
\end{eqnarray*}

$H|_{t=0} = \log tr_{\omega_0} (\omega_T)$ is bounded from above.
By the maximum principle,

$$ H \leq - C\int_0^T \log (T-t) dt + C\leq C' .$$

The proposition is proved as $\varphi$ is uniformly bounded from above.

\qed
\end{proof}


\begin{corollary} There exists $C>0$ such that for $(t,z)\in [0, T)\times X$

\begin{equation}
\omega_t \geq C \pi^* \omega_{FS}.
\end{equation}

\end{corollary}


\begin{theorem} \label{lowerbound} Let $g(t)$ be the K\"ahler metric associated with $\omega(t, \cdot)$. Then there exists $C>0$, such that for $t\in [0, T)$,

\begin{equation}
diam(X, g(t)) \geq C.
\end{equation}

\end{theorem}

\begin{proof}  Since $Y$ is normal, the singular set of $Y$ is an analytic subvariety of $Y$ of codimension greater than one. Let $x_0$ be a point in the nonsingular part $Y_{reg}$ of $Y$. Then there exists a geodesic ball $B_{2r} (x_0, g_{FS}) \subset Y_{reg}$ of radius $r>0$ centered at $z_0$ with respect to the Fubini-Study metric $g_{FS}$ restricted on $Y$. There exist $x_1$ and $x_2 \in B_{r}(x_0, g_{FS})$ such  $$d_{Y, g_{FS}}(x_1, x_2) \geq r ,$$ where $d_{Y, g_{FS}}(x_1, x_2)$ is the distance between $x_1$ and $x_2$ on $Y$ with respect to $g_{FS}$.

Choose $z_1 \in \pi^{-1} (x_1)$ and $z_2 \in \pi^{-1}(x_2)$.
Then for $i = 1$,  $2$,  let $$d_{\pi^{-1}(B_{2r}(x_0, g_{FS})), g(t)}( z_i, \pi^{-1} \partial(  B_{2r} (x_0, g_{FS})  ))$$ be the distance from $z_i$ to $ \pi^{-1}(\partial(  B_{2r} (x_0, g_{FS} ) )$ in $\pi^{-1}(B_{2r}(x_0, g_{FS})$ with respect to $g(t)$ and $d_{Y, g_{FS}} (x_i , \partial B_{2r} (x_0, g_{FS}))$ be the distance from $x_i$ to $\partial B_{2r} (x_0, g_{FS})$ on $Y$ with respect to $g_{FS}$. Then
$$ d_{\pi^{-1}(B_{2r}(x_0, g_{FS})), g(t)}( z_i, \pi^{-1} \partial(  B_{2r} (x_0, g_{FS})  )) \geq C d_{Y, g_{FS}} (x_i , \partial B_{2r} (x_0, g_{FS})) \geq C r,$$
and so
$$d_{X, g(t)} (z_1, z_2) \geq  C d_{Y, g_{FS}}(x_1, x_2) \geq Cr .$$

\qed\end{proof}

\section{The case when $K_X$ is nef}

When $K_X$ is nef, it follows from Proposition \ref{maxtime} that the unnormalized K\"ahler-Ricci flow (\ref{ricciflow}) has long time existence.

\begin{proposition}\label{cy} Let $X$ be $n$-dimensional projective manifold of $\kod(X)=0$.  If $K_X$ is semi-ample, then there exists $C>0$ depending on $g_0$ such that along the K\"ahler-Ricci flow (\ref{ricciflow})

\begin{equation}
diam(X, g(t) ) \geq C.
\end{equation}
\end{proposition}

\begin{proof} If $\kod(X)=0$ and $K_X$ is semi-ample, $c_1(X)=0$. The proposition is a result of Cao \cite{Cao}, \cite{Y1}.

\qed\end{proof}

\begin{proposition} \label{mainprop} Let $X$ be $n$-dimensional projective manifold of $\kod(X)>0$.  If $K_X$ is semi-ample, then there exists $C>0$ depending on $g_0$ such that along the K\"ahler-Ricci flow (\ref{ricciflow})

\begin{equation}
diam(X, g(t) ) \geq C \sqrt{t}.
\end{equation}

\end{proposition}

\begin{proof} We consider the following normalized K\"ahler-Ricci flow when $K_X$ is nef.

\begin{equation}\label{normalflow}
\left\{
\begin{array}{rl}
&{ \displaystyle \frac{ \partial }{\partial s} \tilde{g} =-Ric( \tilde{g} )  - \tilde{g} }\\
&\\
& \tilde{g} |_{s=0}= g_0 .
\end{array} \right.
\end{equation}
The relation between the solution $g(t)$ of the unnormalized K\"ahler-Ricci flow (\ref{ricciflow}) and the solution $\tilde{g}(s)$ of the normalized flow (\ref{normalflow}) is given by
$$t = e^s -1,$$
$$g(t) = e^{s} \tilde{g} (s)= (t+1) \tilde{g}(\log(t+1)).$$

If $X$ is of general type, there exists a birational holomorphic map $f: X \rightarrow X_{can}$ from $X$ to its canonical model $X_{can}$ as $R(X, K_X)$ is finitely generated. Let $E$ be the set of all the points on $X$ such that $f$ is not isomorphic.   Then $E$ is a subvariety of $X$ and it is proved in \cite{Ts}, \cite{TZha} that $g(t)$ converges in $C^{\infty}(X\setminus E)$ as $t \rightarrow \infty$ along the normalized flow (\ref{normalflow}). Since the limiting singular K\"ahler-Einstein metric is smooth and non-degenerate on $X\setminus E$, $diam(X, \tilde{g}(t))$ is then uniformly bounded below from $0$.

If $ 0< k=\kod(X) < n$, then there exists a unique holomorphic map $f: X \rightarrow X_{can}\in \mathbf{CP}^N$, where $X_{can}$ is the canonical model of $X$ of $\dim X_{can}  =k$. Let $$X_{can}^{\circ} = \{ x \in X_{can}~|~ x ~ is ~nonsingular ~and ~f^{-1}(x)~is ~ nonsingular\}.$$ Let $h$ be the pullback of the Fubini-Study metric on $X_{can}$. 	Then it is proved in \cite{SoT1}, \cite{SoT2} by the parabolic Schwarz lemma, that for any $K\subset \subset X_{can}^{\circ}$, there exists $C_K>0$ such that on $[0, \infty)\times K$,
\begin{equation}
\tilde{g}(t) \geq C_K f^* h
\end{equation}
along the normalized K\"ahler-Ricci flow (\ref{normalflow}).
Therefore, $diam(X, \tilde{g}(t))$ is uniformly bounded from below from $0$.

Combining the above estimates, for any initial metric $g_0$, there exists $C>0$ such that along the normalized K\"ahler-Ricci flow (\ref{normalflow}),
\begin{equation}
 diam(X, \tilde{g}(s) ) \geq C.
 \end{equation}
The proposition is proved as $g(t)= (t+1) \tilde{g}(\log(t+1))$.
\qed\end{proof}



\medskip

Now we are able to prove the main theorems.

\bigskip
\noindent{\bf Proof of Theorem \ref{mainth2} }If $K_X$ is not nef and $K^{-1}_X$ is not ample, $L = H + TK_X$ has positive Iitaka dimension at the singular time $t=T$.  If $K_X^{-1}$ is ample and $[\omega_0]$ is not proportional to $K_X^{-1}$, $K_X$ is not nef and $L = H + TK_X$ has positive Iitaka dimension at the singular time $t=T$ by Proposition \ref{fano}. Then the proposition follows by combining Theorem \ref{lowerbound}, Proposition \ref{cy} and  Proposition \ref{mainprop}.  \qed

\bigskip

\noindent{\bf Proof of Theorem \ref{mainth1} }
The K\"ahler-Ricci flow (\ref{ricciflow}) has long time existence if $K_X$ is nef. Theorem \ref{mainth1} follows easily from Theorem \ref{mainth2} and Theorem \ref{perelman}. \qed


\bigskip
\noindent
{\bf Acknowledgements.} \  The author is grateful to Professor D.H. Phong for his advice, encouragement and support.  The author would like to thank:  Professor G. Tian for bringing his article \cite{T2} to his attention as well as his useful suggestion; V. Tosatti for a number of very enlightening discussions; Professor J. Sturm, G. Szekelyhidi, B. Weinkove and Y. Yuan for helpful suggestions.

\bigskip
\bigskip

\noindent $^{*}$ Department of Mathematics \\
Rutgers University, Piscataway, NJ 08854\\

\end{document}